\documentclass[]{svmult}

\usepackage{type1cm}          
\usepackage{makeidx}          
\usepackage{graphicx}         
\usepackage{multicol}         
\usepackage[bottom]{footmisc} 

\usepackage{newtxtext}        %
\usepackage[varvw]{newtxmath} 

\usepackage{algorithm}
\usepackage{algpseudocode}
\usepackage{setspace}
\usepackage{subcaption}
\usepackage[normalem]{ulem}
\usepackage[hidelinks]{hyperref} 

\makeindex

\begin{document}

\title*{Towards the Use of Anderson Acceleration in Coupled Transport-Gyrokinetic Turbulence Simulations}
\author{
David J. Gardner\orcidID{0000-0002-7993-8282},
Lynda L. LoDestro\orcidID{0000-0001-6503-4189},
Carol S. Woodward\orcidID{0000-0002-6502-8659}}
\authorrunning{D. J. Gardner, L. L. LoDestro, and C. S. Woodward}
\institute{
David J. Gardner, Lynda L. LoDestro, and Carol S. Woodward
\at Lawrence Livermore National Laboratory, Livermore, CA, USA\\
\email{{gardner48, lodestro1, woodward6}@llnl.gov}
}

\maketitle

\abstract{Predicting the behavior of a magnetically confined fusion plasma over long time periods requires methods that can bridge the difference between transport and turbulent time scales.
The nonlinear transport  solver, Tango, enables simulations of very long times, in particular to steady state, by advancing each process independently with different time step sizes and couples them through a relaxed iteration scheme.
We examine the use of Anderson Acceleration (AA) to reduce the total number of coupling iterations required by interfacing Tango with the AA implementation, including several  extensions to AA, provided by the KINSOL nonlinear solver package in SUNDIALS.
The ability to easily enable and adjust algorithmic options through KINSOL allows for rapid experimentation to evaluate different approaches with minimal effort.
Additionally, we leverage the GPTune library to automate the optimization of algorithmic parameters within KINSOL.
We show that AA can enable faster convergence in stiff and very stiff tests cases without noise present and in all cases, including with noisy fluxes, increases robustness and reduces sensitivity to the choice of relaxation strength.}

\section{Introduction}
\label{s:intro}

Accurately simulating turbulent transport is a critical aspect of modeling magnetically confined fusion plasmas.
Turbulent fluxes have a strong nonlinear dependence on the gradients of averaged fields and are often the primary mechanism for transporting particles and energy from the plasma interior.
A major challenge in simulations of very long times, in particular steady state, is the significant difference in characteristic time scales for transport (seconds) and turbulence (microseconds) in magnetic-fusion-energy burning plasmas \cite{post2004report}.
This challenge is further complicated by the presence of noise in the averaged transport fluxes arising from turbulent simulations.

Evolving both processes at the same time step is prohibitively expensive; however, as the dynamics on each time scale are well separated, each process can be advanced independently at its own rate and coupled through an iteration scheme to bridge the scale gap.
Turbulence model evaluations are the dominant computational expense, and the speed at which the coupling iteration converges determines the number of turbulence model evaluations.
In this work, we aim to reduce the total number of coupling iterations necessary and thus lower the overall simulation cost by accelerating the convergence of the iterative method.
To this end, we interface the nonlinear transport solver, Tango \cite{parker2018investigation,parker2018bringing}, with the KINSOL nonlinear solver package from SUNDIALS \cite{gardner2022enabling,hindmarsh2005sundials} to examine the use of Anderson Acceleration (AA) and extensions to AA for improving the convergence rate of the coupling iteration.

The remainder of this paper is organized as follows.
In Section \ref{s:coupling}, we overview turbulent transport simulations for magnetically confined fusion plasmas and the coupling method applied.
We introduce the AA algorithm and detail the various extensions we consider to improve performance in Section \ref{s:anderson}.
In Section \ref{s:impl} we overview the Tango code for coupled simulations and its interfacing with the AA implementation from KINSOL as well as the use of GPTune to optimize parameter selection.
In Section \ref{s:results} we discuss results from numerical simulations with Tango in both stiff and very stiff settings along with the impact of noisy fluxes on performance.
Finally, in Section \ref{s:conclustions}, we provide concluding remarks and directions for future investigation.

\section{Coupling Transport \& Turbulence}
\label{s:coupling}

Simulating turbulence-driven transport in toroidal magnetic-fusion-energy plasmas with near-first-principles accuracy requires solving a set of gyrokinetic equations describing the distribution of different species in a 5D phase space (3 physical- plus 2 velocity-space dimensions). 5D turbulence due to fast time-scale instabilities drives transport of averaged fields (i.e., density, momentum, and energy, which are constant on magnetic flux surfaces and thus spatially 1D in the magnetic flux coordinate) across the magnetic surfaces confining the plasma. Because the transport time scales are much slower than the turbulence time scales, a formal separation-of-scales approach is possible, permitting a multiscale computational solution which proceeds by coupling separate equation sets for the fast and slow scales \cite{shestakov2003self}.
For this work, we utilize a simplified model for a single species to evaluate algorithms for turbulent transport without the expense of 5D simulation or the complexities of the toroidal geometry.
Consider the transport equation,
\begin{equation}
\label{eq:transport}
\frac{\partial p}{\partial t} + \frac{\partial q}{\partial x} = S,
\end{equation}
where $p(t,x)$ is the plasma profile, $q(p(t,x))$ is the average turbulent flux, and $S(t,x)$ includes all local sources.
Because the turbulence-driven transport is stiff, an implicit time-advance is required.
Discretizing \eqref{eq:transport} in space with finite differences and in time with backward Euler, we obtain the implicit system,
\begin{equation}
\label{eq:transport_discrete}
p^{(n+1)} + H \mathcal{D}_x q^{(n+1)} = p^{(n)} + H S^{(n+1)},
\end{equation}
for a time step from $t^{(n)}$ to $t^{(n+1)}$ with step size $H$ where superscripts denote the time at which a quantity is approximated, and $\mathcal{D}_x$ is the difference operator for the finite-difference method.

In Tango, we solve the implicit system \eqref{eq:transport_discrete} using the LoDestro method as detailed in \cite{crotinger1997corsica} and \cite{shestakov2003self}.
The method is given by the iteration,
\begin{equation}
\label{eq:transport_iteration}
p_{k+1} - H \mathcal{D}_x( D_k \mathcal{D}_x p_{k+1} - c_k p_{k+1}) = p^{(n)} + H S^{(n+1)},
\end{equation}
where $p_{k+1}$ is the new iterate for $p^{(n+1)}$ and $q^{(n+1)}$ is split into diffusive and convective terms, with coefficients
\begin{equation}
\label{eq:transport_iteration_coeffs}
D_k = -\theta_k q(p_k) / \mathcal{D}_x p_k
\quad\mathrm{and}\quad
c_k = (1 - \theta_k) q(p_k) / p_k,
\end{equation}
respectively.
The parameter $\theta_k$ controls the diffusive-convective split at each mesh point.
The purpose of the split is to keep $D_{k}$ positive (assuring a parabolic linear equation) and within bounds.
It is assumed that the transport physics is diffusive at the boundaries; while it is often diffusive throughout the confined plasma volume, there are important exceptions, which a general treatment must accommodate. With a purely diffusive ($c_k=0$) representation of $q$, profile maxima can cause violent spatial swings in $D_k$ while iterating; also, converged solutions with the flux running uphill (anti-diffusion, possible in sub-volumes with global physics) are unobtainable. While Tango currently employs only the flux-splitting technique from \cite{crotinger1997corsica} and \cite{shestakov2003self} to handle these complications, \cite{crotinger1997corsica} and \cite{shestakov2003self} provide a second technique for representing $q$, which worked about equally well for the problems studied there.
The important point is that however $q$ is represented, it is not a model of the physics (which is unknown, apart from the assumption at the boundaries); yet, assuming the iteration converges, the numerical $q$ supplied from the turbulence code must be recovered, as \eqref{eq:transport_iteration} and \eqref{eq:transport_iteration_coeffs} have been constructed to do here.

The value of $\theta_k$ is chosen using the default strategy in Tango,
\begin{equation}
    \theta_k =
    \begin{cases}
      0   & \text{where}\quad \hat{D}_k < D_{\text{min}} \\
      \frac{1}{2} \left(1 + \frac{D_{\text{max}} - \hat{D}_k}{D_{\text{max}} - D_{\text{min}}} \right) & \text{where}\quad D_{\text{min}} \leq \hat{D}_k \leq D_{\text{max}} \;\text{and}\; |\mathcal{D}_x p_k| < 10 \\
      \frac{1}{2} & \text{where}\quad \hat{D}_k > D_{\text{max}} \;\text{and}\; |\mathcal{D}_x p_k| < 10 \\
      1 & \text{otherwise}
   \end{cases},
\end{equation}
with $\hat{D}_k = -q(p_k) / \mathcal{D}_x p_k$, $D_{\text{min}} = 10^{-5}$, and $D_{\text{max}} = 10^{13}$.
Obtaining the new profile in \eqref{eq:transport_iteration} requires solving the linear system,
\begin{equation}
\label{eq:lodestro_ls}
M_k p_{k+1} = b,
\end{equation}
within each iteration, where $M_k = I - H \mathcal{D}_x( D_k \mathcal{D}_x - c_k )$ and $b = p^{(n)} + H S^{(n+1)}$.
In practice, relaxation of the iteration is required to ensure convergence; i.e., $p_{k+1} = \beta p_{k+1} + (1 - \beta) p_{k}$. The overall method is outlined in Algorithm~\ref{alg:lodestro}.
\begin{algorithm}[H]
    \caption{LoDestro Method }\label{alg:lodestro}
    \hspace*{\algorithmicindent} \textbf{Input:} $p_0$, $\beta$, and tol
    \hspace*{\algorithmicindent} \textbf{Output:} $p^*$
    \begin{algorithmic}[1]
        \For{$k = 0,1,\ldots,k_{\rm{max}}$}
            \State Compute the flux $q(p_k)$ and split it into $D_k$ and $c_k$ \label{alg:g_start}
            \State Construct $M_k$ and solve $M_k p_{k+1} = b$ \label{alg:g_end}
            \State Relax the profile $p_{k+1} = \beta p_{k+1} + (1 - \beta) p_k$
            \State\algorithmicif\ $\| p_{k+1} - p_{k} \| < \text{tol}$ \Return $p_{k+1}$ as $p^*$
        \EndFor
        \end{algorithmic}
\end{algorithm}

\begin{figure}[htb!]
\centering
\includegraphics[width=0.65\textwidth]{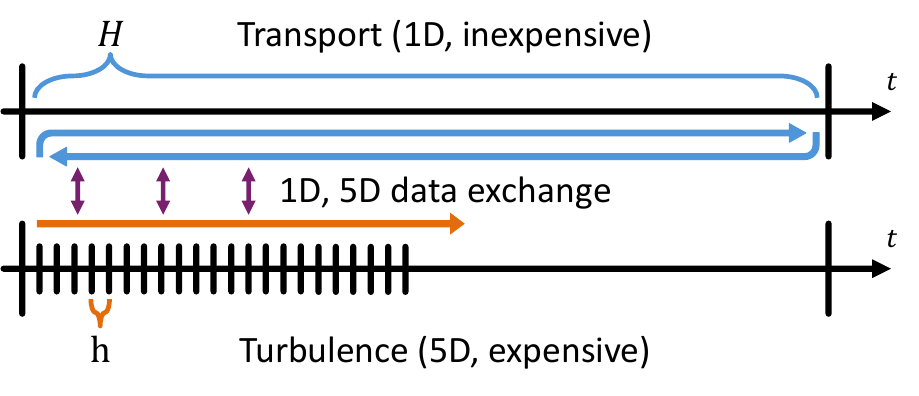}
\caption{
In Algorithm~\ref{alg:lodestro}, the} transport code (1D) uses a large implicit time step that requires iterating to convergence. During this iteration, the turbulence model (5D for physics runs; 1D in this paper) advances with a small explicit time step. The two models periodically exchange fluxes and profiles until the iteration has converged.
\label{fig:lodestro}
\end{figure}

As illustrated in Fig.~\ref{fig:lodestro}, the coupled simulation alternates between the transport iteration and the turbulent evolution. The turbulent fluxes are obtained from a 5D global (entire confined plasma) gyrokinetic code running with an explicit step size, $h \ll H$.
Within each transport iteration the models exchange information. The transport code provides profiles to the turbulence code which, running one long (on the fast time scale) simulation, then continues advancing with updated 1D profiles for a set time, after which it returns fluxes back to the transport code.
The number of fast-scale time steps per data exchange is kept to a multiple of fluctuation periods, rather than the number to converge to a statistical steady state given the latest 1D profile; this provides a large savings in run-time per iteration, mitigating the slow rate of fixed-point convergence.
Evolving the turbulence model requires far more computational resources and expense than solving the linear system \eqref{eq:lodestro_ls} in the transport iteration.
Minimizing the number of iterations for convergence is, then, the critical issue for overall performance.

References \cite{crotinger1997corsica} and \cite{shestakov2003self} include applications of the LoDestro method with a two-dimensional turbulence model demonstrating speed-ups in reaching turbulent steady states.
The first-ever multiscale results coupling to a global 5D gyrokinetic turbulence simulation were obtained by Tango and presented in \cite{parker2018bringing}.
Three-field coupling studies of current experiments are presented in \cite{Di-Siena_2022}, on the ASDEX Upgrade tokamak; and in \cite{Banon-Navarro_2023}, on ion-temperature clamping in the Wendelstein 7-X stellarator.
The first of these includes a section on timings and speed ups.

\section{Anderson Acceleration}
\label{s:anderson}

The LoDestro method can be viewed as a damped fixed-point iteration given by
\begin{equation}
    p_{k+1} = \beta G(p_k) + (1 - \beta) p_k,
\end{equation}
where $G(p)$ comprises steps \ref{alg:g_start} and \ref{alg:g_end} in Algorithm~\ref{alg:lodestro}.
This iteration may converge quite slowly ($\beta \ll 1$ may be required for stability of the iterations), and, thus, we are interested in methods for accelerating convergence.
One such method is Anderson Acceleration (AA) \cite{anderson1965iterative, kelley2018numerical} which computes an iteration update that approximately minimizes the linearized fixed-point residual over a subspace constructed from prior iterations.
In a linear setting and under certain conditions, AA is essentially equivalent to the generalized minimal residual (GMRES) method \cite{walker2011anderson}, and, with a sufficiently good initial iterate, AA performs no worse than a fixed-point iteration \cite{toth2015convergence}.

The AA algorithm with fixed damping, $\beta$, and depth (subspace size), $m$, is given in Algorithm~\ref{alg:anderson}.
Following \cite{walker2011anderson}, the optimization problem in step \ref{alg:aa_simple_minimize} for the weights, $\gamma_k \in \mathbb{R}^{m_k}$, to compute the new iterate is presented as an unconstrained linear least-squares problem.
Utilizing this formulation, rather than the constrained form of the problem, allows for an efficient solution via a QR factorization of $\mathcal{F}_k$ i.e., $R_k \gamma_k = Q^T_k f_k$ where $Q_k \in \mathbb{R}^{N \times m_k}$ is an orthogonal matrix, $R_k \in \mathbb{R}^{m_k \times N}$ is an upper triangular matrix, and $f_k \in \mathbb{R}^{N}$ is the nonlinear residual, $G(p_k) - p_k$.
Updating the factorization and solving the linear system as the residual history changes is typically a minimal cost compared to the function evaluation \cite{lockhart2022performance,loffeld2016considerations} and thus not a major contributor to the overall iteration run time.

\begin{algorithm}[]
    \setstretch{1.1}
    \caption{Anderson Acceleration}\label{alg:anderson}
    \hspace*{\algorithmicindent} \textbf{Input:} $p_0$, $m_{\rm{max}}$, $k_{\rm{max}}$, tol, and $\beta$
    \hspace*{\algorithmicindent} \textbf{Output:} $p^*$
    \begin{algorithmic}[1]
        \State Set $p_1 = \beta G(p_0) + (1 - \beta) p_0$
        \For{$k = 1,\ldots,k_{\rm{max}}$}
            \State Set $m_k = \min\{k, m_{\rm{max}}\}$
            \State Set $\mathcal{G}_k = [\Delta G_{k-m_k}, \ldots, \Delta G_{k-1}]$ where $\Delta G_i = G(p_{i+1}) - G(p_i)$
            \State Set $\mathcal{F}_k = [\Delta f_{k - m_k}, \ldots, \Delta f_{k-1}]$ where $\Delta f_i = f_{i+1} - f_i$ and $f_i = G(p_i) - p_i$
            \State Determine $\gamma_k$ that minimizes $\| f_k - \mathcal{F}_k \gamma_k \|_2$ \label{alg:aa_simple_minimize}
            \State Set $p_{k+1} = G(p_k) - \mathcal{G}_k \gamma_k - (1 - \beta)(f_k - Q_k R_k \gamma_k)$
            \State\algorithmicif\ $\| p_{k+1} - p_{k} \| < \text{tol}$ \Return $p_{k+1}$ as $p^*$
        \EndFor
        \end{algorithmic}
\end{algorithm}

In this work, we consider several extensions to AA, including delaying the start of acceleration as well as non-stationary variants with variable damping and/or depth, leading to Algorithm~\ref{alg:anderson_ext}.
When the initial iterate is far from the fixed-point, delaying acceleration by $d$ iterations can be advantageous in order to obtain a better initial iterate before enabling AA.
Different approaches to adaptive damping in \cite{evans2020proof} and \cite{chen2024non} demonstrated significant improvements in convergence when compared to fixed damping values.
The method in \cite{chen2024non} computes an ``optimal'' damping value by solving a minimization problem at the cost of two additional function evaluations. To avoid extra evaluations, we use the heuristic approach from \cite{evans2020proof} with the damping value
\begin{equation} \label{eq:adapt_beta}
    \beta_k = 0.9 - \omega_{\beta} \Gamma_k
    \quad
    \mathrm{with}
    \quad
    \Gamma_k = \sqrt{1 - (\|Q^Tf_k\|/\|f_k\|)^2},
\end{equation}
where $\Gamma_k$ is the AA convergence gain in the $k$-th iteration, $Q$ is readily available from the QR solve in the AA optimization step, and $\omega_{\beta}$ is a tuning parameter.
In \cite{chen2022composite} the adaptive damping algorithm from \cite{chen2024non} is combined with a composite version of AA as a means of adapting the depth and further improves the convergence rate.
A different adaptivity approach is taken in \cite{pollock2023filtering} where a filtering method is applied to remove columns from the history matrix, $\mathcal{F}_k$, in order to control the condition number of the least-squares problem.
This method also leads to improved convergence. In this work, we utilize a strategy similar to that employed in \cite{pollock2021anderson} where the depth is determined by the magnitude of the current residual with
\begin{equation} \label{eq:adapt_m}
m_k = \min\{\tilde{m}_k, m_{k-1} + 1, m_{max}\}
\quad
\mathrm{with}
\quad
\tilde{m}_k = \max\{ 0, \lfloor -\log_{10} ( \omega_m \| f_k \| ) \rfloor \}
\end{equation}
where $m_0 = 0$ and $\omega_m$ is a turning parameter.

\begin{algorithm}[]
    \setstretch{1.1}
    \caption{Anderson Acceleration with delay and variable damping and depth}\label{alg:anderson_ext}
    \hspace*{\algorithmicindent} \textbf{Input:} $p_0$, $m_{\rm{max}}$, $k_{\rm{max}}$, tol, $d$, and $\beta$
    \hspace*{\algorithmicindent} \textbf{Output:} $p^*$
    \begin{algorithmic}[1]
        \State Set $p_1 = \beta G(p_0) + (1 - \beta) p_0$
        \For{$k = 1,\ldots,k_{\rm{max}}$}
            \If{$k \leq d$}
            \State Set $p_{k+1} = \beta G(p_k) + (1 - \beta) p_k$
            \Else
            \State Select $m_k$ using \eqref{eq:adapt_m}
            \State Set $\mathcal{G}_k = [\Delta G_{k-m_k}, \ldots, \Delta G_{k-1}]$ where $\Delta G_i = G(p_{i+1}) - G(p_i)$
            \State Set $\mathcal{F}_k = [\Delta f_{k - m_k}, \ldots, \Delta f_{k-1}]$ where $\Delta f_i = f_{i+1} - f_i$ and $f_i = G(p_i) - p_i$
            \State Determine $\gamma_k$ that minimizes $\| f_k - \mathcal{F}_k \gamma_k \|_2$ \label{alg:aa_minimize}
            \State Select $\beta_k$ using \eqref{eq:adapt_beta}
            \State Set $p_{k+1} = G(p_k) - \mathcal{G}_k \gamma_k - (1 - \beta_k)(f_k - Q_k R_k\gamma_k)$
            \EndIf
            \State\algorithmicif\ $\| p_{k+1} - p_{k} \| < \text{tol}$ \Return $p_{k+1}$ as $p^*$
        \EndFor
        \end{algorithmic}
\end{algorithm}

\section{Implementation}
\label{s:impl}

To evaluate whether the application of AA improves the convergence rate of the LoDestro method, we utilize the Tango 1D transport solver \cite{parker2018investigation,parker2018bringing} interfaced with the KINSOL nonlinear system solver from the SUNDIALS library \cite{gardner2022enabling,hindmarsh2005sundials}.
Tango is an open-source Python code that both solves the 1D transport equation \eqref{eq:transport_iteration} by way of the linear system \eqref{eq:lodestro_ls} and implements the coupling Algorithm~\ref{alg:lodestro}.
The turbulent flux models available for coupling include the global gyrokinetic code GENE \cite{genecode,goerler2011global}, as well as analytic models for rapid algorithmic development. For this work we utilize one of the latter options, and note it corresponds to running the turbulence code to saturation each coupling exchange.

KINSOL provides a number of algorithms for solving nonlinear algebraic systems including an Anderson-accelerated fixed-point iteration. AA in KINSOL leverages several implementation optimizations for greater performance such as efficient QR factorization updates and solves for computing the acceleration weights \cite{fang2009two, kresse1996efficiency, walker2011anderson} or the use of low synchronization orthogonalization methods within the QR solve to minimize the number of global reductions \cite{lockhart2022performance}.
For this work, we have extended KINSOL's existing support for delay and constant damping with a fixed depth to allow for adapting the depth and/or damping.
These new capabilities will be included in an upcoming SUNDIALS release.
Through the interface between Tango and KINSOL, we can readily access an efficient implementation of AA and easily experiment with different variants of AA with minimal modification to existing Tango code.

The interfacing between Tango and KINSOL primarily consists of creating a class within Tango with a method to evaluate $G(p)$ using Tango's native functions for computing and splitting fluxes and setting up and solving the linear system.
Additional class methods call KINSOL functions to configure the solver based on command line options and solve the nonlinear system.
As SUNDIALS packages are written in C, we leverage SWIG \cite{beazley1996swig} to generate Python interfaces that align with the C API while allowing Tango to access data stored in SUNDIALS vectors as NumPy arrays without having to copy data between KINSOL and Tango data structures.

Algorithm~\ref{alg:anderson_ext} includes a number of parameters that can have a significant impact on performance.
The optimal combination of settings is unclear \textit{a priori} and depends on the application and problem considered.
To automate parameter selection we utilize the GPTune package \cite{liu2021gptune}, a Python-based auto-tuning framework using Bayesian optimization methodologies with support for multi-task learning, multi-objective tuning,  sensitivity analysis, and many other features.
For this work we use the single-task learning autotuning (SLA) capability in GPTune for a single-objective function.
That is, we consider one configuration of the problem \eqref{eq:transport} when optimizing the algorithm parameters to minimize the number of iterations to reach a given solution tolerance.
The problem class created for interfacing with KINSOL also streamlines using GPUTune, as the objective function only needs to construct the Tango problem class given a dictionary of input options, call the solve method, and return the number of solver iterations.
Additionally, the driver program defines the input, output, and parameter spaces, as well as any constraints, to create the tuning problem. Once the problem is defined, the driver creates the GPTune object and runs the SLA method which calls Tango as needed through the objective function.

\section{Numerical Results}
\label{s:results}

In this section we present results applying AA to the LoDestro method for a nonlinear diffusion problem from \cite{shestakov2003self} that mimics the challenges observed in turbulence-driven transport simulations.
The problem is given by \eqref{eq:transport} on the interval $x \in [0,1]$ with the analytic flux
\begin{equation} \label{eq:shestakov_flux}
    q = -D \frac{\partial p}{\partial x} + \epsilon \quad \text{where} \quad D = \left( \frac{1}{p} \frac{\partial p}{\partial x} \right)^{r}.
\end{equation}

As briefly mentioned in previous sections,
turbulent transport in tokamaks is often due to instabilities driven by gradients of the time-averaged 1D temperature and other profiles; and these instabilities often have thresholds---critical values for their onset. The resulting transport then drives the system to flatter gradients, which pushes plasma operation toward marginal stability, where the growth rates rise extremely rapidly from zero. This situation leads to very stiff transport as operation near the critical gradients is approached from above \cite{Garbet_2004}.
We model this high level of stiffness by increasing the value of $r$ in \eqref{eq:shestakov_flux}.
The system is stiff with $r = 2$ and we use $r = 10$ to represent very stiff tests relevant to magnetic-fusion-energy problems.

The 5D turbulence code illustrated in Fig.~\ref{fig:lodestro} reaches a saturated state with short spatial- and temporal-scale fluctuations. It is the electric fields involved in these fluctuations that kick particles across the magnetic surfaces, i.e., cause the transport.
Even when a heat or particle flux is integrated over four dimensions, the 1D flux still has short-scale structure.
We include the effects of these fluctuations on our coupling problem with the ``noise" term, $\epsilon$, in \eqref{eq:shestakov_flux}.
We use the same spatially correlated, temporally white Gaussian noise model studied in \cite{parker2018investigation} to mimic the fluctuations observed in GENE simulations (see Figs.~4a and~5b in \cite{parker2018investigation} for a comparison of the model with fluxes from GENE).
Below, we first test cases without noise ($\epsilon = 0$) before considering a problem with noisy fluxes.

The source term in \eqref{eq:transport} is defined as
\begin{equation}
    S(x) =
    \begin{cases}
      1 & \text{for}\quad 0 \leq x \leq 0.1 \\
      0 & \text{for}\quad 0.1 < x \leq 1
   \end{cases},
\end{equation}
and the boundary conditions are $p'(0) = 0$ and $p(1) = 0.01$.
For the LoDestro method, the relaxation strength depends on the degree of stiffness, and the default value is $\beta = 0.6/r$.
AA setups with zero depth ($m = 0$) and constant damping correspond to applying the unaccelerated LoDestro method as implemented in Tango and serves as the baseline comparison point.

\subsection{A Stiff Test Case}
\label{s:stiff_test}

We first consider the stiff case ($r = 2$) and the impact of different AA configurations on the convergence of the iteration.

\subsubsection{Fixed Damping, Depth, and Delay}
\label{s:stiff_case_fixed_options}

\begin{figure*}[htb!]
    \centering
    \begin{subfigure}[t]{0.5\textwidth}
        \centering
        \includegraphics[width=0.98\textwidth]{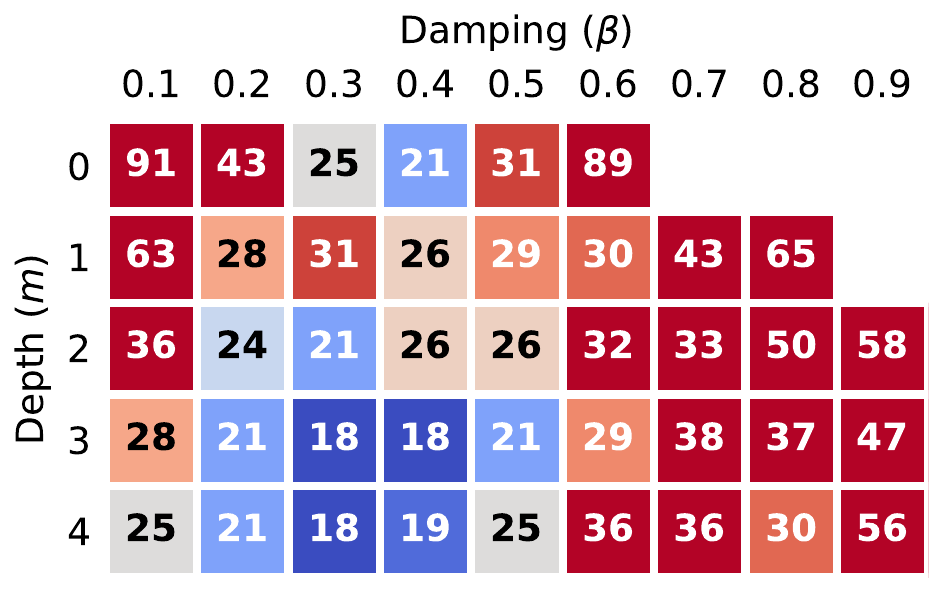}
        \caption{Without Delay}
        \label{fig:heatmap_no_delay}
    \end{subfigure}%
    \begin{subfigure}[t]{0.5\textwidth}
        \centering
        \includegraphics[width=0.98\textwidth]{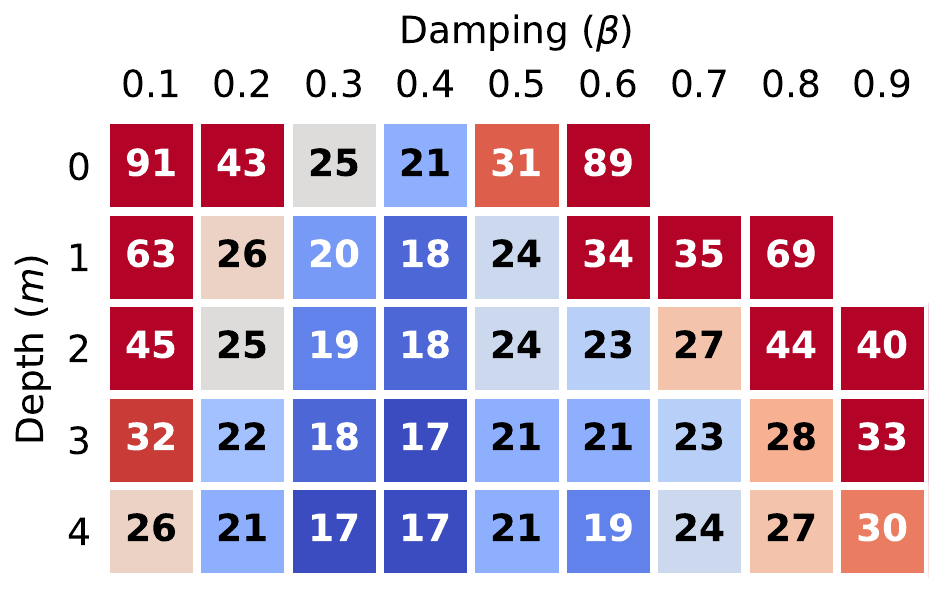}
        \caption{Delay ($d$) AA by 1 iteration}
        \label{fig:heatmap_delay}
    \end{subfigure}
    \caption{Heat map of the number of iterations to reach a residual of $10^{-11}$. Empty cells correspond to failed solves. Without a delay, $m \geq 2$ is needed to see a benefit in iteration counts from AA. With a one iteration delay, AA consistently reduces the number of iterations and decreases sensitivity to the value of $\beta$.}
    \label{fig:heatmap}
\end{figure*}

To begin, we consider AA with static damping, $\beta$, and depth, $m$, sweeping over the parameter space to better understand how each impacts performance.
Fig.~\ref{fig:heatmap} shows a heat map for the number of iterations to reach a residual of $10^{-11}$ for different $\beta$ and $m$ values with and without a delay, $d$, before starting AA.
From the first row ($m = 0$, no acceleration), we see the default, $\beta = 0.3$, is near optimal in this case, but $\beta = 0.4$ can reduce the number of iterations by $16\%$.
However, increasing $\beta$ too much leads to unphysical values or divergence (both denoted by empty cells).

Without delay, adding acceleration, $m = 1$, gives some robustness in the choice of $\beta$, enabling convergence for larger values than possible without AA. However, the convergence rate is improved only for $\beta$ far from the optimal value.
Increasing the depth, $m > 1$, yields more consistent gains with AA.
Comparing the best results, AA takes approximately $14\%$ fewer iterations.
For a non-optimal $\beta$, larger reductions in iteration count are realized due to the greater robustness and reduced sensitivity to the choice of $\beta$.
Improved robustness in $\beta$ with AA has also been observed in other multiphysics coupling applications \cite{hamilton2016assessment}.

As the first iterate diverges from final solution before steadily approaching the correct value at each subsequent iteration, introducing a small delay before starting AA significantly strengthens the performance in all aspects. Generally, the results continue to improve with increased depth but with diminishing returns, especially near the optimal $\beta$. Delaying AA further can improve the results slightly but also reduces robustness with larger $\beta$ and may increase iteration counts with smaller $\beta$. The size of the delay necessary will depend on the quality of the initial iterate, and the best depth value is often problem dependent.

\subsubsection{Optimized Fixed Damping, Depth, and Delay}
\label{s:stiff_case_optimize_fixed_options}

\begin{figure}[htb!]
\sidecaption
\includegraphics[width=7.5cm]{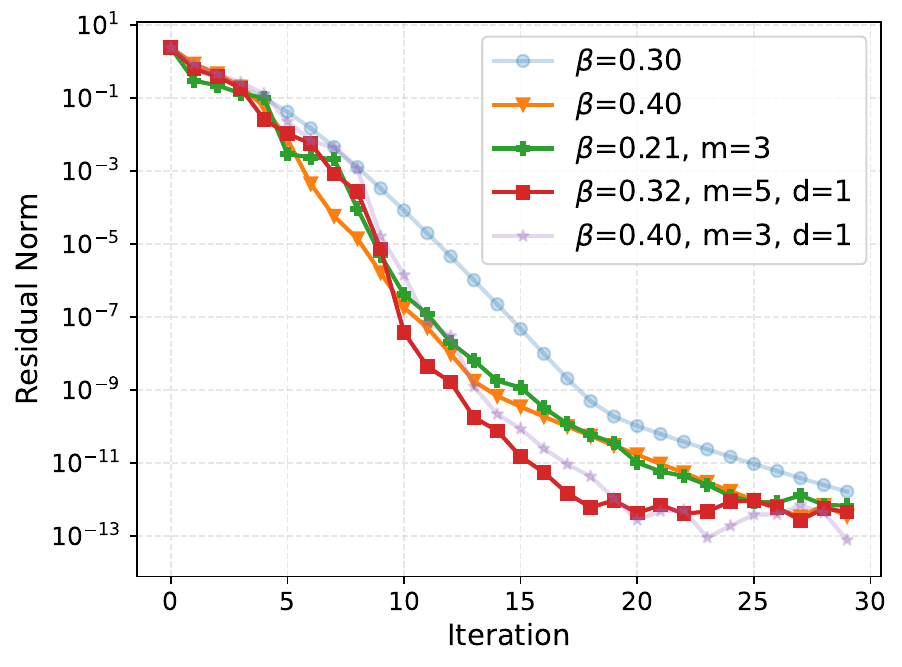}
\caption{Residual history for autotuned configurations. The default (faded circles) and one setup from Fig.~\ref{fig:heatmap_delay} (faded stars) are included for reference. The autotuned results align well with the parameter sweep in Fig.~\ref{fig:heatmap}.
\break \rule[1.5cm]{0pt}{2pt}
}
\label{fig:stiff_case_optimized}
\end{figure}

In the prior section we manually explored the parameter space to understand the impact of different options on algorithmic performance.
In this section we investigate optimizing these selections with GPTune and examine the convergence history.
We consider three sets of parameters to optimize: $\beta$ with $m=0$; $\beta$ and $m$; and $\beta$, $m$, and $d$.
Each set is optimized for three residual tolerances ($10^{-6}$, $10^{-8}$, and $10^{-11}$). We set the initial number of GPTune parameter-space samples to 10 and the total number of samples to 50.
The convergence histories for the best parameter values from each set are shown in Fig.~\ref{fig:stiff_case_optimized}.

In the first few iterations there is little difference between the various options with smaller $\beta$ values performing slightly better.
At more moderate residuals, optimizing $\beta$ alone is one or two iterations better than the AA results.
However, for values below $10^{-7}$ the AA setups give the best results with the optimized parameter set edging out the best result from Fig.~\ref{fig:heatmap_delay} by one or two iterations.
The benefit of AA increases for values less than $10^{-9}$ where the convergence of the unaccelerated version slows.  This result indicates that AA delivers accelerated convergence to a tight tolerance.
Overall the optimized parameter sets from GPTune closely align with the best setups found with the manual sweep.

\subsubsection{Adaptive Damping and Depth}
\label{s:stiff_case_adapt_options}

\begin{figure}[htb!]
\sidecaption
\includegraphics[width=7.5cm]{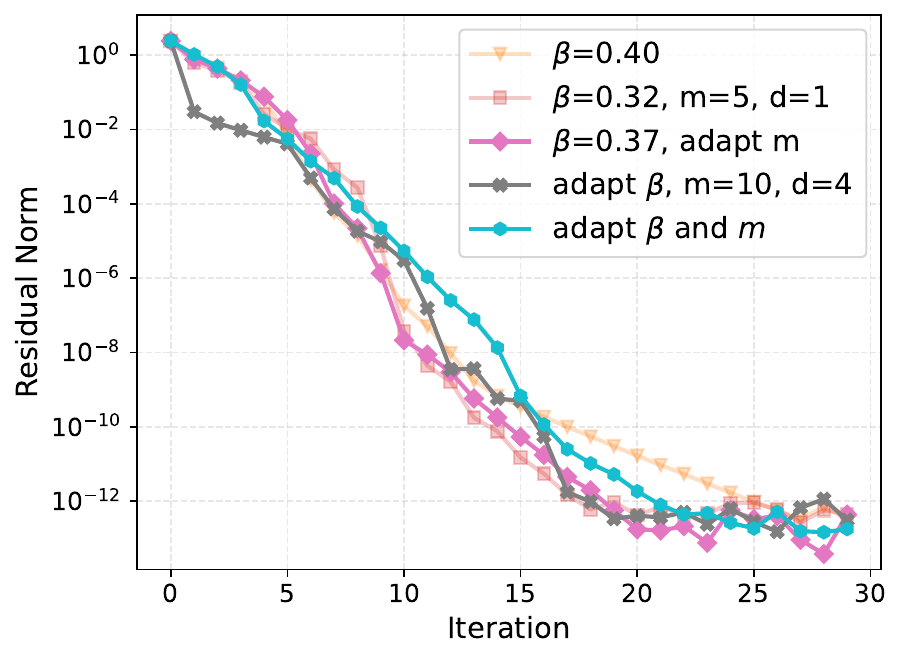}
\caption{Residual history for autotuned configurations with adaptive $\beta$ and/or $m$.
The adaptive results closely follow the autotuned fixed parameter results in Fig.~\ref{fig:stiff_case_optimized} (faded triangles and squares included for reference).
\break \rule[1.5cm]{0pt}{2pt}
}
\label{fig:stiff_case_adaptive}
\end{figure}

\begin{figure*}[htb!]
    \centering
    \begin{subfigure}[t]{0.5\textwidth}
        \centering
        \includegraphics[width=0.95\textwidth]{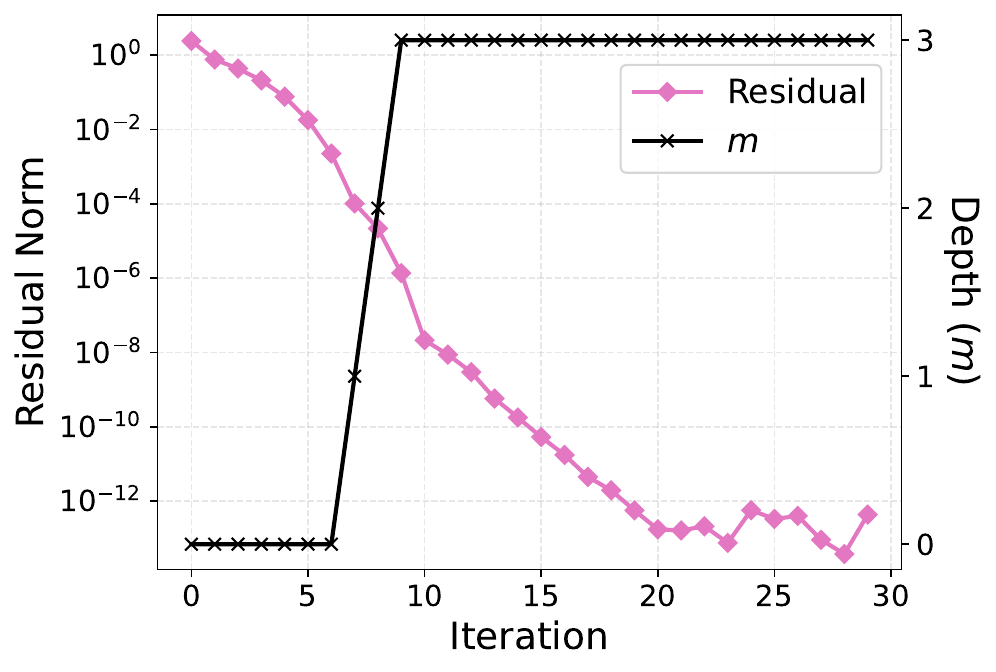}
        \caption{$\beta = 0.37$ and adaptive $m$}
        \label{fig:adaptive_depth_history}
    \end{subfigure}%
    \begin{subfigure}[t]{0.5\textwidth}
        \centering
        \includegraphics[width=0.95\textwidth]{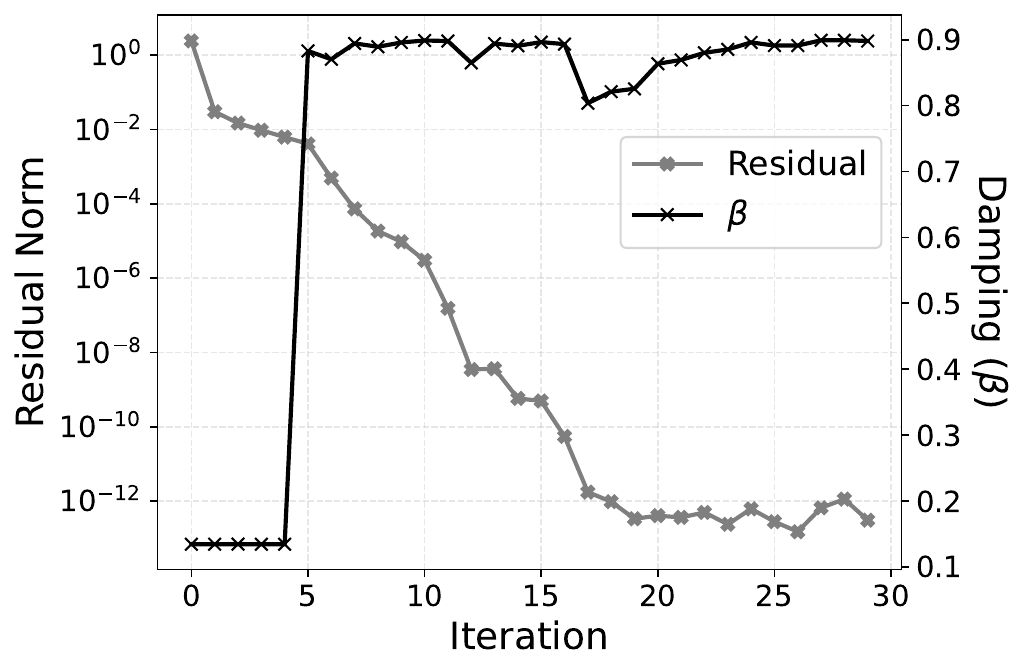}
        \caption{Adaptive $\beta$, $m=10$, and $d=4$}
        \label{fig:adaptive_damping_history}
    \end{subfigure}%
    \par\medskip 
    \begin{subfigure}[t]{0.5\textwidth}
        \centering
        \includegraphics[width=0.95\textwidth]{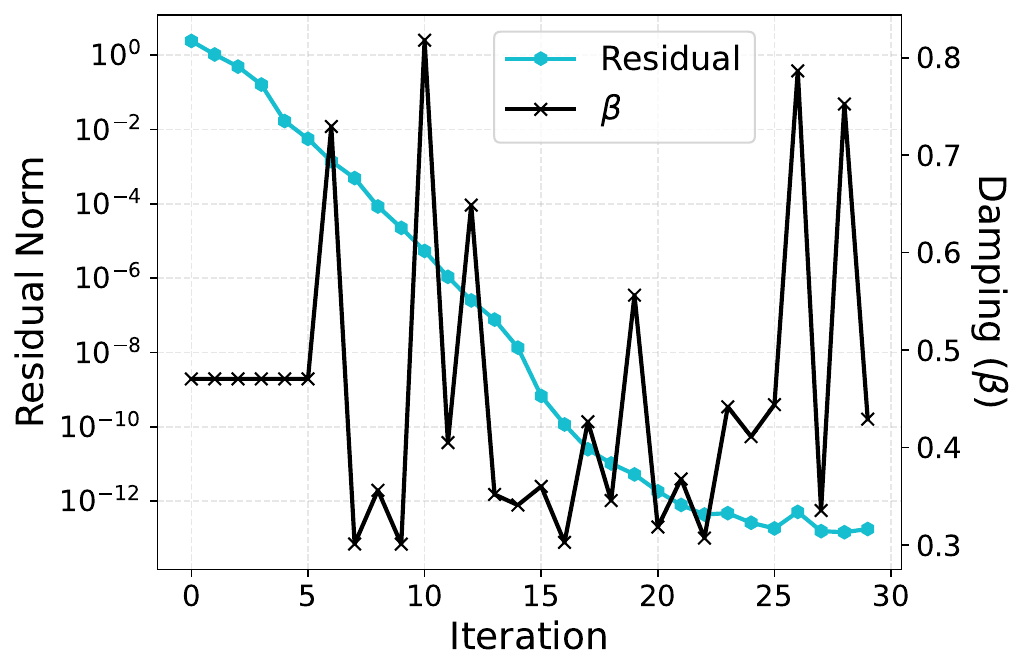}
        \caption{Adaptive $\beta$ and $m$}
        \label{fig:adaptive_damping_and_depth_history_beta}
    \end{subfigure}%
    \begin{subfigure}[t]{0.5\textwidth}
        \centering
        \includegraphics[width=0.95\textwidth]{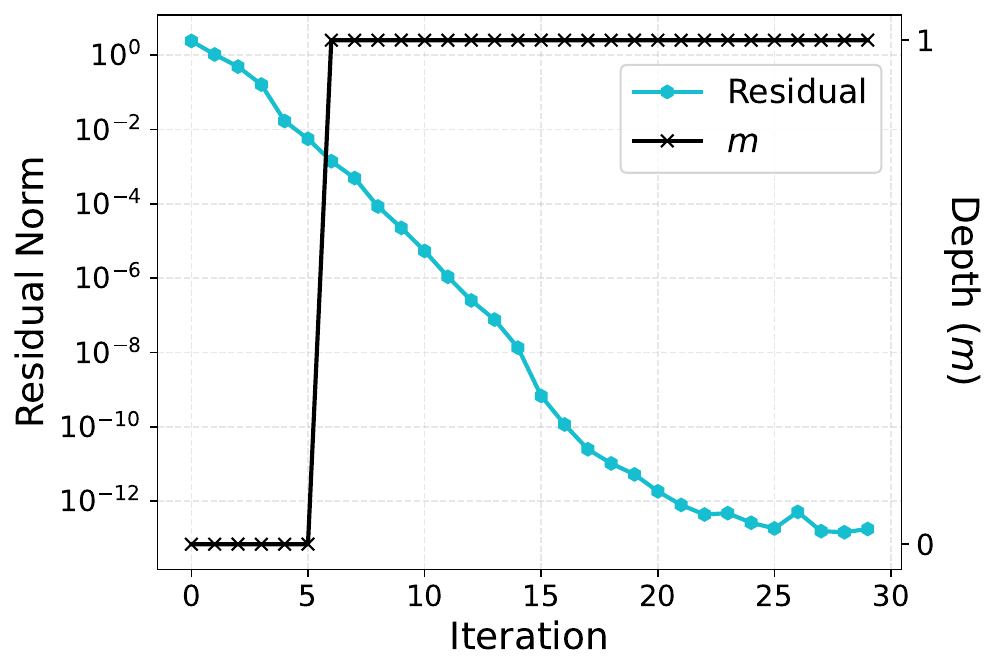}
        \caption{Adaptive $\beta$ and $m$}
        \label{fig:adaptive_damping_and_depth_history_m}
    \end{subfigure}
    \caption{Residual (left axis) and parameter history (right axis) for adaptive $\beta$ and/or $m$ configurations. Adapting $m$ automatically selects the AA delay, and adapting $\beta$ alone leads to values near the max allowed. Adapting $\beta$ and $m$ together leads to large, frequent changes in $\beta$ slightly degrading performance compared to other setups.}
    \label{fig:adaptive_history}
\end{figure*}

In the previous section we saw that the best parameter configurations depend on the desired residual norm.
Performance also depends significantly on whether the initial iterate is in the asymptotic convergence regime of the iterative scheme.
In this section we evaluate AA with adaptive damping and/or depth to find setups that perform well across residual regimes.
We again consider three sets of parameters to optimize: $\beta$ with adaptive $m$; $m$ and $d$ with adaptive $\beta$, and adapting $\beta$ and $m$ together.
We only optimize $d$ in the fixed-$m$ case as the strategy for adapting $m$ automatically finds an optimum $d$ by initially selecting $m = 0$.
When adapting $\beta$ and/or $m$, GPTune optimizes the factors $\omega_\beta$ and $\omega_m$ in \eqref{eq:adapt_beta} and \eqref{eq:adapt_m}, respectively.
As before, each set is optimized for the same three residual tolerances and we use the same numbers of initial and total parameter-space samples.
The convergence history for the best parameter values from each set are shown in Fig.~\ref{fig:stiff_case_adaptive}, and history of $\beta$ or $m$ from the adaptive cases are given in Fig.~\ref{fig:adaptive_history}.

Optimizing a fixed $\beta$ and an adaptive $m$ gives similar results to optimizing fixed $\beta$, $m$, and $d$.
Fig.~\ref{fig:adaptive_depth_history} shows that adaptive $m$ delays AA six iterations and has a max $m$ of three.
Given the low sensitivity observed in Fig.~\ref{fig:heatmap_delay}, it is not surprising that, although the $m$ and $d$ are slightly different, the convergence history only differs by one iteration.
Optimizing an adaptive $\beta$ with fixed $m$ and $d$ has improved convergence at larger residuals, and, after five iterations, the setup is within two iterations of the best results.
Fig.~\ref{fig:adaptive_damping_history} shows that the performance gain at larger residuals is due to using a small $\beta$ value initially before immediately growing to almost the maximum $\beta$ once AA starts.
This behavior occurs because an optimum fixed $\beta$ is computed by GPTune in the delay region before AA starts, as $\Gamma_k$ is not available to estimate $\beta$.
Once AA starts, the iteration switches to adaptive values and keeps $\beta$ close to the max allowed.

Optimizing both $\beta$ and $m$, is the least performant of the three configurations but, due to the low sensitivity to $\beta$ and $m$, the results are still relatively close (within five iterations) of the other cases.
Again adapting $m$ automatically delays AA (five iterations) before settling on a max depth.
Unlike when adapting $\beta$ alone, there are large oscillations in the $\beta$ after the initial delay.
This variation in $\beta$ appears to be the cause of lower overall performance and suggests that considering the gain history to reduce large, frequent changes in $\beta$ may improve results.

Overall, in this stiff test case there are modest gains in numbers of iterations with AA compared to optimized damping values at large to moderate residual values.
For residuals below $10^{-8}$ AA has a larger impact as convergence of the iteration without AA slows.
However, AA greatly increases robustness and decreases sensitivity to the choice of $\beta$ with a small delay in AA increasing these effects.

\subsection{A Very Stiff Test Case}
\label{s:very_stiff_test}

\begin{figure}[htb!]
\sidecaption
\includegraphics[width=7.5cm]{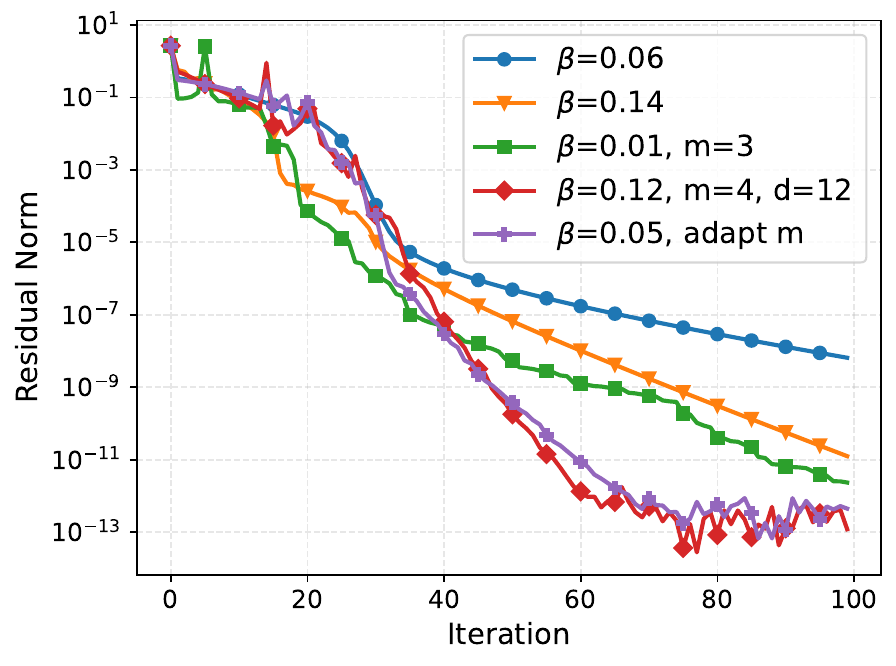}
\caption{Residual history for autotuned setups with the default iteration ($\beta=0.06$), optimized damping ($\beta=0.14$), and various AA configurations. Markers are placed every 5 iterations. For a wide range of residuals, AA leads to reduced numbers of iterations, especially for smaller residual values.
\break \rule[1.0cm]{0pt}{2pt}
}
\label{fig:very_stiff_test_optimized}
\end{figure}

We now consider a very stiff case ($r = 10$) relevant to magnetic-fusion-energy problems and examine the impact of AA on convergence with optimized fixed and adaptive parameter configurations as in the prior two sections.
Fig.~\ref{fig:very_stiff_test_optimized} shows the convergence history for the default iteration ($\beta = 0.06$), an optimized damping value ($\beta = 0.14$), and three of the best performing AA configurations.

As in the stiff case, optimizing $\beta$ gives better results than the default $\beta$ and, in this case, is generally the best option for larger residuals.
Enabling AA with $\beta = 0.01$ and $m = 3$ consistently outperforms optimizing $\beta$ alone for residuals below $10^{-4}$.
Again, delaying AA generally leads to better results overall and, given the larger initial residual and stronger damping applied in this case, the best AA configurations tend to have a longer delay.
For residuals below $10^{-7}$, configurations with fixed $\beta$, $m$, and $d$ or fixed $\beta$ and adaptive $m$ are the most efficient with results that differ by less than six iterations.
With an adaptive $m$, AA is automatically delayed by twelve iterations, the same as in the fixed-delay case, but allows more depth going up to $m=7$.
Unlike in the stiff test, the depth does not immediately ramp up after the delay but instead keeps $m$ constant for 2 to 7 iterations before increasing $m$.
The stair-step behavior of $m$ reflects the overall slower convergence in this case, which requires additional iterations before a residual reduction triggers an increase in $m$ using the heuristic \eqref{eq:adapt_m}.
Adapting $\beta$ can also improve results but, like in the stiff case, does not  generally perform as well as other setups due to large oscillations in $\beta$.
Overall the results with this very stiff test are improved when comparing against the stiff case, with AA producing more significant gains over a larger range of residuals.

\subsection{Adding Noise in a Stiff Test Case}
\label{s:noise}

Finally, we revisit the stiff test ($r = 2$) when accounting for fluctuations in $q$ that mimic the noise observed with GENE i.e., $\epsilon \neq 0$ in \eqref{eq:shestakov_flux}.
Fig.~\ref{fig:noisy_flux} shows the residual history for the default setup ($\beta = 0.3$) and a more optimal configuration ($\beta = 0.4$), with and without AA, along with a heat map for other combinations of $\beta$ and $m$.
\newsavebox\mybox
\savebox{\mybox}{\includegraphics[width=0.49\textwidth]{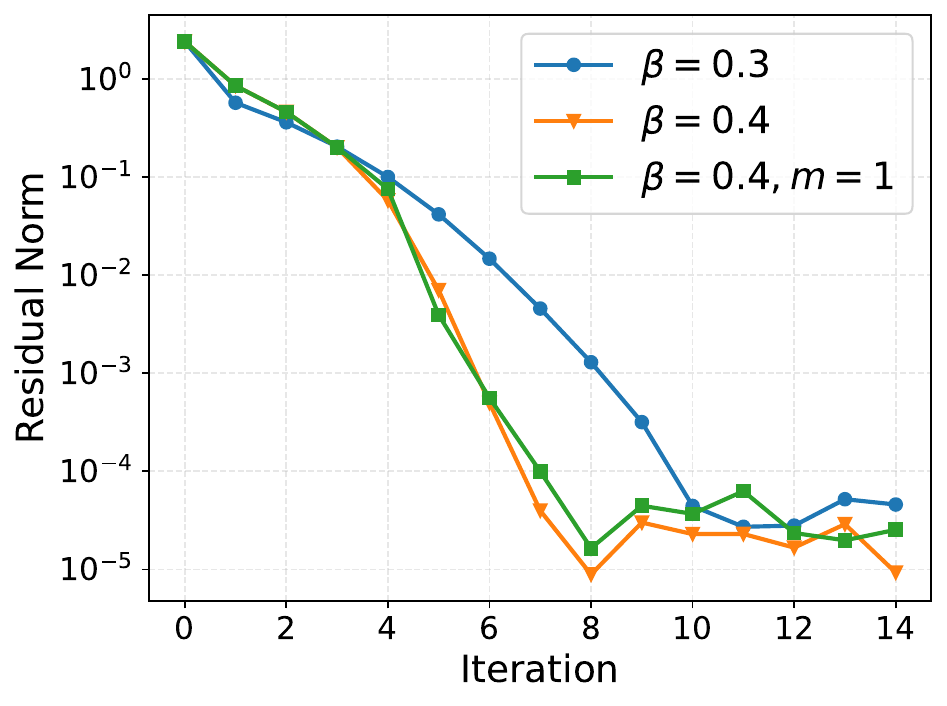}}
\begin{figure}
    \centering
    \subcaptionbox*{}[.49\textwidth]{\usebox{\mybox}}\hfill%
    \subcaptionbox*{}[.49\textwidth]{%
    \vbox to \ht\mybox{%
    \includegraphics[width=0.49\textwidth]{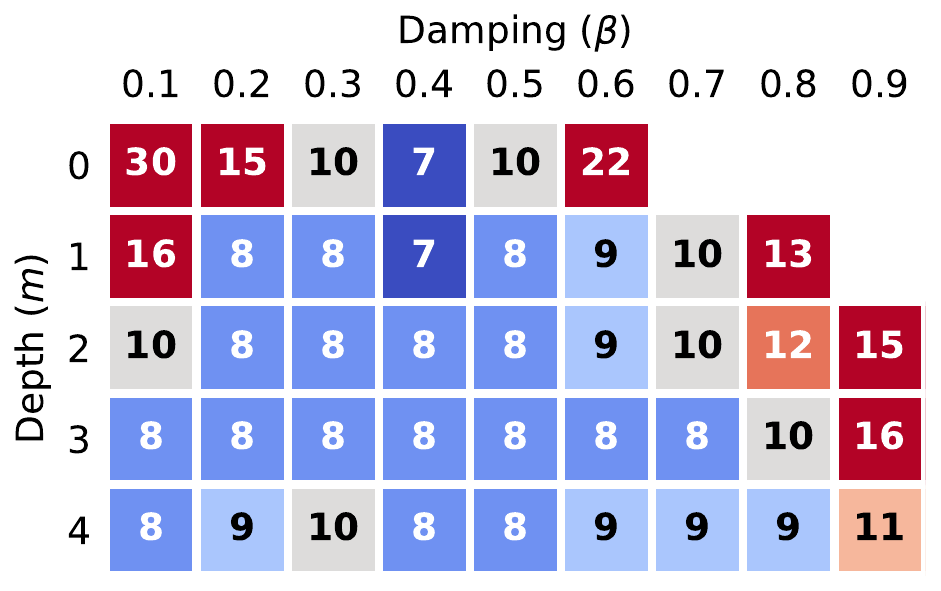}}}
    \vspace{-1.5\baselineskip}
    \caption{Residual history and heat map for the number of iterations to reach a residual of $10^{-4}$ both with a delay of two iterations. AA increases robustness and reduces sensitivity to $\beta$ but only lowers iteration counts away from the optimal $\beta$.}
    \label{fig:noisy_flux}
\end{figure}
Unlike in the cases without noise, the performance with AA is similar to that without AA when using the best $\beta$ value.
This lack of improvement is due to the higher, noise-dependent residual floor at which the iteration stagnates as noted in \cite{parker2018investigation} for the LoDestro method and similarly observed in \cite{toth2017local} when applying AA to problems with inaccurate function evaluations.
While the residual floor limits the benefits in lowering iteration counts, AA still improves overall robustness and reduces sensitivity to $\beta$, which leads to faster convergence for non-optimal $\beta$ values.
This reduced sensitivity is significant as it shows that with AA, less manual tuning of $\beta$ is necessary for an efficient solution approach.

\section{Conclusions}
\label{s:conclustions}

In this work we used the KINSOL nonlinear solver library in SUNDIALS to evaluate the potential of several variations of Anderson acceleration to improve the convergence of the LoDestro method for coupled turbulence and transport simulations of a magnetically confined fusion plasma.
Combining the flexibility of KINSOL with the GPTune autotuning library enabled rapid evaluation of the parameter space to compare the impact of different options on algorithmic performance.

For all test cases, AA allows for convergence with larger relaxation parameter values and decreases sensitivity to the choice of the relaxation parameter.
Comparing to the optimal damping value, AA can provide a small decrease in the number of iterations for moderate residual values with the benefit increasing for smaller target residuals leading to a reduction of 20\% or more in the number of iterations.
Moreover, these benefits improve as the problem difficulty increases, going from stiff to very stiff tests; and for non-optimal damping values, the improvements are larger.
Utilizing adaptive depth with AA provides a means for automatically selecting the delay and maximum allowed depth with results generally matching the best setups with optimized fixed parameter values.
We note that use of GPTune removes the need for time-consuming hand optimizations of the parameters.
When noisy fluxes are considered, AA does not provide an improvement in iterations with optimal damping due to the higher residual floor \cite{parker2018investigation}.
However, AA still increases robustness and reduces sensitivity to the damping parameter, which lowers iteration counts for non-optimal damping.
Improvements in the adaptivity strategies for damping and depth (or the alternative methods explored in \cite{chen2024non}, \cite{evans2020proof}, and \cite{pollock2023filtering}) could increase the benefits of AA at larger residual values and will be considered in later studies.

In future efforts we plan to assess the performance of AA with a modified analytic flux model that reflects the time dependence of the turbulence calculation as it is coupled to the transport code (see Fig.~\ref{fig:lodestro}).
We will also apply recently developed approaches for analyzing AA convergence \cite{de2022linear} to better understand how the presence of noise and the addition of time dependence in the coupling impact convergence of the algorithm.
Additionally, we will conduct studies to quantify parameter sensitivity in AA using GPUTune.
Finally, the interfaces to KINSOL will allow for testing new variants of AA as they are integrated into the library, such as alternating and composite AA \cite{chen2022composite}, without modifying Tango.


\begin{acknowledgement}
The authors wish to acknowledge the initial explorations applying Anderson acceleration in Tango by Jeff Parker, Jack G. Paulson, and H. Hunter Schwartz and thank them for their help in using Tango and interfacing with KINSOL. We would also like to thank Cody J. Balos for his work in creating Python interfaces to KINSOL.

Support for this work was provided in part by the U.S. Department of Energy, Office of Science, Office of Advanced Scientific Computing Research, Scientific Discovery through Advanced Computing (SciDAC) Program through the Frameworks, Algorithms, and Scalable Technologies for Mathematics (FASTMath) Institute, under Lawrence Livermore National Laboratory subcontract B626484 and DOE award DE-SC0021354.

This material is based upon work supported by the U.S. Department of Energy, Office of Science, Office of Advanced Scientific Computing Research and Office of Fusion Energy Sciences, Scientific Discovery through Advanced Computing (SciDAC) program.

This work was performed under the auspices of the U.S. Department of Energy by Lawrence Livermore National Laboratory under Contract DE-AC52-07NA27344. LLNL-PROC-864392.

\end{acknowledgement}
\ethics{Competing Interests}{
The authors have no conflicts of interest to declare that are relevant to the content of this chapter.}

\bibliographystyle{spmpsci}
\bibliography{references}

\end{document}